\documentclass[11pt,a4paper]{article}

\usepackage[margin=27mm]{geometry}
\usepackage{amsmath,amssymb,amsthm,mathtools,mathrsfs}
\usepackage{microtype}
\usepackage{needspace}
\usepackage[hidelinks,unicode]{hyperref}
\usepackage{aliascnt}
\usepackage[nameinlink,noabbrev]{cleveref}
\usepackage{indentfirst}
\usepackage{enumitem}

\allowdisplaybreaks[2]

\theoremstyle{plain}
\newtheorem{theorem}{Theorem}[section]
\newaliascnt{lemma}{theorem}
\newtheorem{lemma}[lemma]{Lemma}
\aliascntresetthe{lemma}
\newaliascnt{proposition}{theorem}
\newtheorem{proposition}[proposition]{Proposition}
\aliascntresetthe{proposition}

\theoremstyle{definition}
\newaliascnt{definition}{theorem}
\newtheorem{definition}[definition]{Definition}
\aliascntresetthe{definition}
\newaliascnt{problem}{theorem}
\newtheorem{problem}[problem]{Problem}
\aliascntresetthe{problem}

\crefname{theorem}{Theorem}{Theorems}
\crefname{lemma}{Lemma}{Lemmas}
\crefname{proposition}{Proposition}{Propositions}
\crefname{definition}{Definition}{Definitions}
\crefname{problem}{Problem}{Problems}
\crefname{section}{Section}{Sections}
\crefname{equation}{Equation}{Equations}

\numberwithin{equation}{section}

\newcommand{\NN}{\mathbb N}
\newcommand{\M}{\mathcal M}
\newcommand{\dGH}{d_{\mathrm{GH}}}
\newcommand{\eps}{\varepsilon}
\DeclareMathOperator{\diam}{diam}
\DeclareMathOperator{\dis}{dis}
\newcommand{\doi}[1]{\href{https://doi.org/#1}{\nolinkurl{doi:#1}}}
\newcommand{\arxiv}[1]{\href{https://arxiv.org/abs/#1}{\nolinkurl{arXiv:#1}}}

\title{Bounded Metric Spaces Embeddable into the Gromov--Hausdorff Space}
\author{Ryo Fukano}
\date{November 30, 2026}

\begin{document}
\maketitle
\begin{NoHyper}
\begingroup
\renewcommand{\thefootnote}{}
\endgroup
\end{NoHyper}

\begin{abstract}
Let $(\M,\dGH)$ be the Gromov--Hausdorff space of isometry classes of nonempty compact metric spaces. 
For separable bounded metric spaces, we give three equivalent intrinsic characterizations of isometric embeddability into $\M$. 
The first is a finite-coloring condition for graphs determined by the metric. 
The second is the closed-ball condition $\alpha_X(\overline B_X(o,r))\le r$, where $\alpha_X$ denotes the Kuratowski measure of noncompactness. 
The third is the equality $\alpha_X(A)=\chi_X(A)$ on every nonempty subset, where $\chi_X$ denotes the Hausdorff measure of noncompactness. 
Thus a bounded metric space is isometrically embeddable into $\M$ if and only if it is separable and satisfies any one of these three conditions.
\end{abstract}

\section{Introduction}\label{sec:intro}

Let $\M$ denote the space of isometry classes of nonempty compact metric spaces equipped with the Gromov--Hausdorff distance. The purpose of this paper is to determine which bounded metric spaces are isometrically embeddable into $\M$.

In a previous paper \cite[Theorem~1.1]{Fukano}, the author proved that every nonempty compact metric space admits an isometric embedding into $\M$. The next problem, posed in \cite[Problem~1.3]{Fukano}, is to characterize the metric spaces that are isometrically embeddable into $\M$. This paper solves this problem for bounded metric spaces and gives intrinsic metric conditions characterizing embeddability.

For the Gromov--Hausdorff class in which noncompact representatives are allowed, every bounded metric space is known to be isometrically embeddable \cite[Theorem~3]{ITclass}. More generally, the same conclusion holds for every extended metric space, in which the distance is allowed to take the value $\infty$ \cite[Theorem~1.1]{FukanoExtended}. The target considered in the present paper is instead the space $\M$, whose representatives are required to be compact metric spaces.

Since $\M$ is separable, every metric space that is isometrically embeddable into $\M$ is separable. We therefore state the bounded characterization for separable bounded metric spaces. The main theorem gives three equivalent intrinsic descriptions of embeddability. They come from finite graph coloring, coverings of closed balls, and equality of two measures of noncompactness.

The closed ball $\overline B_X(o,r)$, condition $(F)$, and the Kuratowski and Hausdorff measures of noncompactness $\alpha_X$ and $\chi_X$ appearing in the theorem below are defined in \cref{sec:prelim}. Our main result is the following intrinsic characterization.

\Needspace{18\baselineskip}
\begin{theorem}\label{thm:main}
Let $X$ be a separable bounded metric space. The following are equivalent.
\begin{enumerate}[label=\textup{(\roman*)},nosep]
\item $X$ is isometrically embeddable into $(\M,\dGH)$.
\item $X$ satisfies condition $(F)$.
\item
\[
 \alpha_X(\overline B_X(o,r))\le r
\]
for every $o\in X$ and $r>0$.
\item
\[
 \alpha_X(A)=\chi_X(A)
\]
for every nonempty subset $A\subseteq X$.
\end{enumerate}
\end{theorem}

The necessary definitions and preliminary facts are collected in \cref{sec:prelim}, and \cref{thm:main} is proved in \cref{sec:proof}.

\section{Preliminaries}\label{sec:prelim}

Throughout the paper, all metric spaces are assumed to be nonempty, and all distances are finite-valued. We write $\NN=\{1,2,\ldots\}$, and sequences are allowed to have repeated terms unless stated otherwise.

\begin{definition}\label{def:metric-notation}
Let $(Z,d_Z)$ be a metric space, $z\in Z$, and $r\ge0$. The closed ball of radius $r$ centered at $z$ is denoted by $\overline B_Z(z,r)=\{w\in Z:d_Z(z,w)\le r\}$. For a nonempty subset $A\subseteq Z$, write $d_Z(z,A)=\inf_{a\in A}d_Z(z,a)$ and $\diam A=\sup_{a,b\in A}d_Z(a,b)$. The subscript is omitted when the ambient metric is clear.
\end{definition}

\begin{definition}\label{def:embedding-separable}
A map $F:(Z,d_Z)\to(W,d_W)$ is an \emph{isometric embedding} if it preserves all distances, that is, if $d_W(F(z),F(z'))=d_Z(z,z')$ for all $z,z'\in Z$. A metric space $Z$ is \emph{separable} if it contains a countable dense subset.
\end{definition}

\begin{definition}\label{def:net-uniform-tb}
Let $Z$ be a metric space and $r>0$. A nonempty finite set $S\subseteq Z$ is an \emph{$r$-net} of $Z$ if $\sup_{z\in Z}d_Z(z,S)\le r$. Define
\[
 N_Z(r)=\min\bigl\{|S|:\varnothing\ne S\subseteq Z\text{ is a finite $r$-net}\bigr\},
\]
with the convention $N_Z(r)=\infty$ if no finite $r$-net exists. The space $Z$ is \emph{totally bounded} if $N_Z(r)<\infty$ for every $r>0$.

A family $\mathscr K$ of compact metric spaces is \emph{uniformly totally bounded} if $\sup_{K\in\mathscr K}\diam K<\infty$ and, for every $r>0$, $\sup_{K\in\mathscr K}N_K(r)<\infty$.
\end{definition}

\begin{definition}\label{def:F}
Let $X$ be a metric space, $o\in X$, and $\eta>0$. Define $\Gamma_{o,\eta}(X)$ to be the graph with vertex set $X$ in which distinct points $x,y$ are adjacent precisely when
\[
 d(x,y)>\max\{d(o,x),d(o,y)\}+\eta.
\]
A \emph{finite coloring} of a graph is a map from its vertex set to a finite set such that adjacent vertices always receive different values. We say that $X$ satisfies \emph{condition $(F)$} if $\Gamma_{o,\eta}(X)$ admits a finite coloring for every $o\in X$ and $\eta>0$.
\end{definition}

\begin{definition}\label{def:mnc}
Let $A$ be a nonempty bounded subset of a metric space $X$. The \emph{Kuratowski measure of noncompactness} $\alpha_X(A)$ is the infimum of all $r\ge0$ such that $A$ can be covered by finitely many nonempty subsets of $X$ of diameter at most $r$. The \emph{Hausdorff measure of noncompactness} $\chi_X(A)$ is the infimum of all $r\ge0$ such that $A$ can be covered by finitely many closed balls of radius $r$ with centers in $X$. No completeness assumption on $X$ is made in either definition.
\end{definition}

If $U\subseteq X$ is nonempty and $\diam U\le r$, then $U\subseteq\overline B_X(u,r)$ for every $u\in U$, while every closed ball of radius $r$ has diameter at most $2r$. Hence
$\chi_X(A)\le\alpha_X(A)\le2\chi_X(A)$ for every nonempty bounded subset $A\subseteq X$.

\begin{definition}\label{def:partitions}
A \emph{finite partition} of a set $S$ is a finite family of pairwise disjoint nonempty subsets whose union is $S$. A partition $\mathcal Q$ \emph{refines} a partition $\mathcal P$ if every $B\in\mathcal Q$ is contained in some $A\in\mathcal P$. The \emph{common refinement} of finitely many partitions of the same set is the partition consisting of all nonempty intersections obtained by choosing one member from each partition.
\end{definition}

\begin{definition}\label{def:gh}
For nonempty compact subsets $A,B$ of a metric space $Z$, their Hausdorff distance is
\[
 d_H^Z(A,B)=\max\left\{\sup_{a\in A}d_Z(a,B),\ \sup_{b\in B}d_Z(b,A)\right\}.
\]
For nonempty compact metric spaces $K,L$, their Gromov--Hausdorff distance is
\[
 \dGH(K,L)=\inf_{(Z,\varphi,\psi)}d_H^Z\bigl(\varphi(K),\psi(L)\bigr),
\]
where the infimum is taken over all triples consisting of a metric space $Z$ and isometric embeddings $\varphi:K\to Z$ and $\psi:L\to Z$. Let $\M$ denote the set of isometry classes of nonempty compact metric spaces. In what follows, a nonempty compact metric space $K$ and its isometry class $[K]\in\M$ are denoted by the same symbol.
\end{definition}

\begin{definition}\label{def:correspondence}
A \emph{correspondence} between sets $K$ and $L$ is a relation $R\subseteq K\times L$ whose two coordinate projections are both surjective. The set of all correspondences between $K$ and $L$ is denoted by $\operatorname{Corr}(K,L)$. For metric spaces $K,L$ and a correspondence $R$ between them, the \emph{distortion} of $R$ is
\[
 \dis R=\sup_{(x,y),(x',y')\in R}|d_K(x,x')-d_L(y,y')|.
\]
\end{definition}

\begin{proposition}[{\cite[Theorem~7.3.25]{BBI}}]\label{prop:gh-correspondence}
For nonempty compact metric spaces $K,L$,
\begin{equation}\label{eq:gh-corr}
 \dGH(K,L)=\frac12\inf_{R\in\operatorname{Corr}(K,L)}\dis R.
\end{equation}
\end{proposition}

\begin{theorem}[Gromov compactness theorem {\cite[Theorem~7.4.15]{BBI}}]\label{thm:gromov-compactness}
Every uniformly totally bounded family of compact metric spaces is relatively compact in $(\M,\dGH)$. Consequently, every sequence in a uniformly totally bounded family of compact metric spaces has a subsequence that converges in the Gromov--Hausdorff distance to a compact metric space.
\end{theorem}

\begin{proposition}\label{prop:gh-basic}
The Gromov--Hausdorff distance is a complete separable metric on $\M$. Moreover, for compact metric spaces $K,L$ and $c>0$,
$|\diam K-\diam L|\le2\dGH(K,L)$ and $\dGH(cK,cL)=c\dGH(K,L)$, where $cK$ denotes the space obtained by multiplying the metric of $K$ by $c$.
\end{proposition}

The fact that $\dGH$ is a metric is \cite[Theorem~7.3.30]{BBI}. Completeness and separability are \cite[Theorems~7.22 and~7.23]{Tuzhilin}. The diameter estimate and scaling identity follow from \cref{prop:gh-correspondence}.

\begin{definition}\label{def:colored-graphs}
A simple graph $H$ consists of a vertex set $V(H)$ and an edge set $E(H)$ of two-element subsets of $V(H)$. We write $st\in E(H)$ for $\{s,t\}\in E(H)$. Two vertices are equivalent if they are equal or can be joined by a finite edge path. The subgraph consisting of one equivalence class and all edges between vertices in that class is called a \emph{connected component}. A \emph{properly 3-colored finite graph} is a finite simple graph $H$ together with a map $\kappa_H:V(H)\to\{0,1,2\}$ taking different values at the endpoints of every edge. For properly 3-colored finite graphs $H,H'$, a \emph{color-preserving graph homomorphism} $f:H\to H'$ is a map $f:V(H)\to V(H')$ that sends edges to edges and satisfies $\kappa_{H'}\circ f=\kappa_H$. The set of all color-preserving graph homomorphisms $H\to H'$ is denoted by $\operatorname{Hom}_{\mathrm{col}}(H,H')$.
\end{definition}

\section{Proof of the Main Theorem}\label{sec:proof}

We first prove the equivalence of the three intrinsic conditions.

\begin{proposition}\label{prop:intrinsic}
Let $X$ be a bounded metric space. The following are equivalent.
\begin{enumerate}[label=\textup{(\alph*)},nosep]
\item $X$ satisfies condition $(F)$.
\item $\alpha_X(\overline B_X(o,r))\le r$ for every $o\in X$ and $r>0$.
\item $\alpha_X(A)=\chi_X(A)$ for every nonempty subset $A\subseteq X$.
\end{enumerate}
\end{proposition}

\begin{proof}
Assume \textup{(a)}. Fix $o\in X$, $r>0$, and $\eta>0$, and choose a finite coloring $c$ of $\Gamma_{o,\eta}(X)$. If $x,y$ lie in the intersection of the same color class with $\overline B_X(o,r)$, then
$d(x,y)\le\max\{d(o,x),d(o,y)\}+\eta\le r+\eta$. Hence the nonempty intersections of the color classes with $\overline B_X(o,r)$ form a finite cover of $\overline B_X(o,r)$ by sets of diameter at most $r+\eta$. Thus $\alpha_X(\overline B_X(o,r))\le r+\eta$, and letting $\eta\downarrow0$ gives \textup{(b)}.

Assume \textup{(b)}, and let $\varnothing\ne A\subseteq X$ and $s>\chi_X(A)$. By the definition of the infimum, there exist $x_1,\ldots,x_N\in X$ such that $A\subseteq\bigcup_{j=1}^N\overline B_X(x_j,s)$. For every $\eps>0$, condition \textup{(b)} implies that each $\overline B_X(x_j,s)$ can be covered by finitely many sets of diameter at most $s+\eps$. Combining these finite covers gives $\alpha_X(A)\le s+\eps$. Letting first $\eps\downarrow0$ and then $s\downarrow\chi_X(A)$ yields $\alpha_X(A)\le\chi_X(A)$. Together with $\chi_X(A)\le\alpha_X(A)$, this proves \textup{(c)}. Conversely, \textup{(c)} gives
$\alpha_X(\overline B_X(o,r))=\chi_X(\overline B_X(o,r))\le r$, so \textup{(b)} follows.

It remains to prove \textup{(b)}$\Rightarrow$\textup{(a)}. Fix $o\in X$ and $\eta>0$, and set $\delta=\eta/3$. Since $X$ is bounded, only finitely many of the sets
$E_j=\{x\in X:j\delta\le d(o,x)<(j+1)\delta\}$, $j=0,1,\ldots$, are nonempty. For every nonempty $E_j$, by \textup{(b)}, $\alpha_X(E_j)\le(j+1)\delta$. Hence $E_j$ admits a finite cover $U_{j,1},\ldots,U_{j,N_j}$ by sets of diameter at most $(j+2)\delta$. Assign to $x\in E_j$ the color
$(j,\min\{\ell:x\in U_{j,\ell}\})$. The coloring has finite range. If $x,y\in E_j$ have the same color, then
\[
 d(x,y)\le(j+2)\delta
 \le\max\{d(o,x),d(o,y)\}+2\delta
 <\max\{d(o,x),d(o,y)\}+\eta.
\]
Thus $x$ and $y$ are not adjacent in $\Gamma_{o,\eta}(X)$, so $\Gamma_{o,\eta}(X)$ is finitely colorable.
\end{proof}

We next prove \cref{thm:main}\textup{(i)}$\Rightarrow$\textup{(iii)}.

\begin{proposition}\label{prop:embedding-to-ball}
If a bounded metric space $X$ is isometrically embeddable into $(\M,\dGH)$, then $X$ is separable and
\[
 \alpha_X(\overline B_X(o,r))\le r
\]
for every $o\in X$ and $r>0$.
\end{proposition}

\begin{proof}
Assume that there is an isometric embedding $F:X\to\M$. Since $\M$ is separable, so is $X$. Fix $o\in X$, $r>0$, and $\delta>0$. Choose a compact representative $K$ of $F(o)$ and cover $K$ by nonempty sets $K_1,\ldots,K_N\subseteq K$ of diameter at most $\delta$. For each $x\in\overline B_X(o,r)$, choose a representative $Y_x$ of $F(x)$. Since $\dGH(K,Y_x)=d(o,x)\le r$, \eqref{eq:gh-corr} gives a correspondence $R_x\in\operatorname{Corr}(K,Y_x)$ satisfying $\dis R_x<2r+\delta$.

For $1\le j\le N$, set
$A_j^x=\{y\in Y_x:\exists z\in K_j,\ (z,y)\in R_x\}$. Surjectivity of the two coordinate projections of $R_x$ implies that every $A_j^x$ is nonempty and that $Y_x=\bigcup_{j=1}^N A_j^x$. The distortion bound gives
$\diam A_j^x\le\diam K_j+\dis R_x<2r+2\delta$ and
$\diam Y_x\le\diam K+\dis R_x<\diam K+2r+\delta$.

Define
\[
 \tau_x(y)=\left(\delta\left\lfloor\frac{d_{Y_x}(y,A_j^x)}\delta\right\rfloor\right)_{j=1}^N,
 \qquad
 \mathcal T_x=\{\tau_x(y):y\in Y_x\}.
\]
All vectors $\tau_x(y)$ belong to the common finite set
$(\delta\mathbb Z)^N\cap[0,\diam K+2r+\delta]^N$, and hence the map $x\mapsto\mathcal T_x$ has finite range. For $v=(v_1,\ldots,v_N)\in\mathbb R^N$, write $\|v\|_\infty=\max_{1\le j\le N}|v_j|$.

For arbitrary $y,z\in Y_x$, distance to a set is $1$-Lipschitz, and the truncation error in each coordinate is nonnegative and less than $\delta$. Hence the difference of the two truncation errors has absolute value less than $\delta$. Therefore
$\|\tau_x(y)-\tau_x(z)\|_\infty\le d_{Y_x}(y,z)+\delta$. Choose $j$ with $z\in A_j^x$. The $j$th coordinate of $\tau_x(z)$ is $0$. Taking the infimum over $A_j^x$ in the triangle inequality gives
\[
\begin{aligned}
 d_{Y_x}(y,z)
 &\le d_{Y_x}(y,A_j^x)+\diam A_j^x\\
 &<\delta\left\lfloor d_{Y_x}(y,A_j^x)/\delta\right\rfloor+2r+3\delta\\
 &\le\|\tau_x(y)-\tau_x(z)\|_\infty+2r+3\delta.
\end{aligned}
\]
Consequently,
\[
 \|\tau_x(y)-\tau_x(z)\|_\infty-\delta
 \le d_{Y_x}(y,z)
 <\|\tau_x(y)-\tau_x(z)\|_\infty+2r+3\delta.
\]

Suppose that $\mathcal T_x=\mathcal T_{x'}$. Then
$\Lambda=\{(y,y')\in Y_x\times Y_{x'}:\tau_x(y)=\tau_{x'}(y')\}$ has surjective coordinate projections and is therefore a correspondence between $Y_x$ and $Y_{x'}$. For any $(y,y'),(z,z')\in\Lambda$, the identities $\tau_x(y)=\tau_{x'}(y')$ and $\tau_x(z)=\tau_{x'}(z')$ imply
\[
\begin{aligned}
 \|\tau_x(y)-\tau_x(z)\|_\infty-\delta
 &\le d_{Y_x}(y,z)
 <\|\tau_x(y)-\tau_x(z)\|_\infty+2r+3\delta,\\
 \|\tau_x(y)-\tau_x(z)\|_\infty-\delta
 &\le d_{Y_{x'}}(y',z')
 <\|\tau_x(y)-\tau_x(z)\|_\infty+2r+3\delta.
\end{aligned}
\]
Hence $\dis\Lambda\le2r+4\delta$, and \eqref{eq:gh-corr} yields
$d(x,x')=\dGH(Y_x,Y_{x'})\le r+2\delta$. The nonempty fibers of the finite-range map $x\mapsto\mathcal T_x$ therefore form a finite cover of $\overline B_X(o,r)$ by sets of diameter at most $r+2\delta$. Thus $\alpha_X(\overline B_X(o,r))\le r+2\delta$, and letting $\delta\downarrow0$ gives $\alpha_X(\overline B_X(o,r))\le r$.
\end{proof}

To prove the remaining implication \cref{thm:main}\textup{(ii)}$\Rightarrow$\textup{(i)}, we prepare four lemmas. The first uses uniform total boundedness and a diagonal argument to pass from finite metric realizations of the distances in a dense sequence to an isometric embedding.

\begin{lemma}\label{lem:finite-approx}
Let $X$ be a separable metric space, and let $(x_i)_{i\in\NN}$ be a dense sequence in $X$. Suppose that finite metric spaces $A_i^{(n)}$ are given for all $n\in\NN$ and $1\le i\le n$, and satisfy
\begin{equation}\label{eq:finite-exact}
 \dGH(A_i^{(n)},A_j^{(n)})=d_X(x_i,x_j)
 \qquad(1\le i,j\le n),
\end{equation}
and that, for each fixed $i$, the family $(A_i^{(n)})_{n\ge i}$ is uniformly totally bounded. Then $X$ is isometrically embeddable into $\M$.
\end{lemma}

\begin{proof}
By \cref{thm:gromov-compactness}, a convergent subsequence can be extracted from $(A_i^{(n)})_{n\ge i}$ for every fixed $i$. Successively passing to nested subsequences for $i=1,2,\ldots$ and then taking a diagonal subsequence $n_k\to\infty$, we may arrange that, for every fixed $i$, the tail with $n_k\ge i$ satisfies $A_i^{(n_k)}\to K_i$ for some compact metric space $K_i$. Passing to the limit in \eqref{eq:finite-exact} gives $\dGH(K_i,K_j)=d_X(x_i,x_j)$. If $x_i=x_j$, then $K_i$ and $K_j$ are isometric, so $x_i\mapsto K_i$ defines a well-defined isometric embedding of the dense subset $\{x_i:i\in\NN\}$.

For $x\in X$, choose a sequence $x_{i_k}\to x$. Then $(K_{i_k})$ is Cauchy in $\M$. By completeness in \cref{prop:gh-basic}, let $F(x)$ be its limit. If $x_{i_k}\to x$ and $x_{j_k}\to x$, then $\dGH(K_{i_k},K_{j_k})=d_X(x_{i_k},x_{j_k})\to0$, so $F(x)$ is independent of the choice. If $x_{i_k}\to x$ and $x_{j_k}\to y$, passing to the limit in the distance identity gives $\dGH(F(x),F(y))=d_X(x,y)$.
\end{proof}

The next lemma constructs the finite partitions needed later from condition $(F)$.

\begin{lemma}\label{lem:partitions}
Let $Y=\{x_i:i\in\NN\}$ be a countable metric space of diameter at most $1$, and suppose that $Y$ satisfies condition $(F)$. For $x\in Y$ and a nonempty set $A\subseteq Y$, define
\begin{equation}\label{eq:H-def}
 H_x(A)=\sup_{\substack{u,v\in A\\z\in Y}}
 \bigl(d(z,u)-d(z,v)-d(x,z)\bigr).
\end{equation}
Then the following properties hold.
\begin{enumerate}[label=\textup{(\alph*)},nosep]
\item $0\le H_x(A)\le\diam A\le1$, and $A\subseteq B$ implies $H_x(A)\le H_x(B)$.
\item $|H_x(A)-H_y(A)|\le d(x,y)$ for all $x,y\in Y$.
\item $\sup_{u,v\in A}|d(x,u)-d(x,v)|\le H_x(A)$.
\item There exists a sequence $(\mathcal P_n)_{n\ge0}$ of finite partitions of $Y$ such that $\mathcal P_0=\{Y\}$, $\mathcal P_n$ refines $\mathcal P_{n-1}$ for $n\ge1$, $\{x_j\}\in\mathcal P_n$ for $1\le j\le n$, and
\begin{equation}\label{eq:partition-decay}
 \max_{A\in\mathcal P_n}H_{x_i}(A)\le2^{-n}
 \qquad(1\le i\le n).
\end{equation}
\end{enumerate}
\end{lemma}

\begin{proof}
Setting $z=x$ and $u=v$ in \eqref{eq:H-def} gives $H_x(A)\ge0$, while the reverse triangle inequality gives $H_x(A)\le\diam A$. Monotonicity follows from inclusion of the ranges over which the supremum is taken. Each term in \eqref{eq:H-def} is $1$-Lipschitz as a function of $x$, proving (b). Taking $z=x$ and considering both orders of $u,v$ gives (c).

To prove (d), first fix $x\in Y$ and $\eta>0$. We construct a finite partition all of whose members $C$ satisfy $H_x(C)\le\eta$. Put $\eps=\eta/3$, and choose a finite coloring $c_0$ of $\Gamma_{x,\eps}(Y)$. Since $d(x,u)\le1$, the map
$u\mapsto(c_0(u),\lfloor d(x,u)/\eps\rfloor)$ has finite range. Let $\mathcal V$ be the family of its nonempty fibers. If $a,b,z\in A\in\mathcal V$, then the values $d(x,u)$ for $u\in A$ differ by less than $\eps$, and hence
$d(a,b)\le\max\{d(x,a),d(x,b)\}+\eps<d(x,z)+2\eps$. Thus $\diam A\le d(x,z)+2\eps$.

Since $0\le d(u,A)\le1$, the quantization map
$u\mapsto(\lfloor d(u,A)/\eps\rfloor)_{A\in\mathcal V}$ has finite range. Consider the finite partition by its nonempty fibers. If $u,v$ lie in the same fiber $C$ and $z\in Y$, choose $A\in\mathcal V$ with $z\in A$. Then $|d(u,A)-d(v,A)|<\eps$. Moreover, $d(z,u)\le d(u,A)+\diam A$ and $d(z,v)\ge d(v,A)$, so
\[
 d(z,u)-d(z,v)-d(x,z)
 \le d(u,A)-d(v,A)+\diam A-d(x,z)<3\eps=\eta.
\]
Taking the supremum over $u,v\in C$ and $z\in Y$ gives $H_x(C)\le\eta$.

Set $\mathcal P_0=\{Y\}$. Suppose $n\ge1$ and $\mathcal P_{n-1}$ has been defined. For each $1\le i\le n$, use the preceding construction to choose a finite partition $\mathcal C_i$ such that every $B\in\mathcal C_i$ satisfies $H_{x_i}(B)\le2^{-n}$. Let $\mathcal P_n$ be the common refinement of $\mathcal P_{n-1}$, $\mathcal C_1,\ldots,\mathcal C_n$, and the partitions $\{\{x_j\},Y\setminus\{x_j\}\}\setminus\{\varnothing\}$ for $1\le j\le n$. Then $\mathcal P_n$ is a finite partition refining $\mathcal P_{n-1}$ and contains $\{x_j\}$ for $1\le j\le n$. For every $A\in\mathcal P_n$ and $1\le i\le n$, some $B\in\mathcal C_i$ contains $A$, so monotonicity in (a) yields $H_{x_i}(A)\le H_{x_i}(B)\le2^{-n}$.
\end{proof}

We next encode inclusion between partition elements by finite graph homomorphisms.

\begin{lemma}\label{lem:poset-graphs}
For every countable partially ordered set $(P,\le)$, there exists a family $(H_p)_{p\in P}$ of properly 3-colored finite graphs, each using all three colors, such that
\begin{equation}\label{eq:hom-order}
 \operatorname{Hom}_{\mathrm{col}}(H_p,H_q)\ne\varnothing
 \quad\Longleftrightarrow\quad p\le q.
\end{equation}
\end{lemma}

\begin{proof}
The case $P=\varnothing$ is trivial, so assume $P\ne\varnothing$. For a positive integer $m$, let $C_m$ be the cycle of length $3m$ with vertex set $(\mathbb Z/m\mathbb Z)\times\{0,1,2\}$ and edges
$(r,0)(r,1)$, $(r,1)(r,2)$, and $(r,2)(r+1,0)$ for $r\in\mathbb Z/m\mathbb Z$, where $(r,c)$ has color $c$. A color-preserving graph homomorphism $C_m\to C_n$ exists if and only if $n\mid m$. Indeed, orient the edges as $0\to1$, $1\to2$, and $2\to0$. A color-preserving graph homomorphism preserves these orientations, and every vertex has a unique outgoing oriented edge. Hence the image of one circuit around $C_m$ returns to its starting point only if $3n\mid3m$. Conversely, if $n\mid m$, then $(r,c)\mapsto(r\bmod n,c)$ is a color-preserving graph homomorphism.

Enumerate $P=\{p_i:i\in I\}$ without repetition, where $I=\{1,\ldots,|P|\}$ if $P$ is finite and $I=\NN$ otherwise. Choose distinct primes $(\pi_i)_{i\in I}$ and set
\[
 N_i=\prod_{\substack{j\le i\\p_i\le p_j}}\pi_j.
\]
Then $N_j\mid N_i$ if and only if $j\le i$ and $p_i\le p_j$. Indeed, if $N_j\mid N_i$, then $\pi_j\mid N_i$, giving $j\le i$ and $p_i\le p_j$. Conversely, if $j\le i$ and $p_i\le p_j$, then every $k$ with $\pi_k\mid N_j$ satisfies $k\le j\le i$ and $p_i\le p_j\le p_k$, so $\pi_k\mid N_i$.

For each $i\in I$, define
\[
 V(H_{p_i})=
 \{(k,r,c):k\le i,\ p_k\le p_i,\ r\in\mathbb Z/N_k\mathbb Z,\ c\in\{0,1,2\}\}.
\]
Two vertices $(k,r,c)$ and $(k,r',c')$ with the same first coordinate are adjacent exactly when $(r,c)(r',c')$ is an edge of $C_{N_k}$. Vertices with different first coordinates are not adjacent. Give $(k,r,c)$ color $c$. Thus each connected component of $H_{p_i}$ is a cycle $C_{N_k}$ with fixed first coordinate $k$, and the component with $k=i$ is present, so all three colors occur.

Suppose $p_i\le p_j$. Define $f:V(H_{p_i})\to V(H_{p_j})$ by
\[
 f(k,r,c)=
 \begin{cases}
 (k,r,c),&k\le j,\\
 (j,r\bmod N_j,c),&k>j.
 \end{cases}
\]
If $k\le j$, then $p_k\le p_i\le p_j$, so the image lies in the component of $H_{p_j}$ with first coordinate $k$. If $k>j$, then $p_k\le p_i\le p_j$ and $j<k$, so $N_j\mid N_k$, and the image lies in the component of $H_{p_j}$ with first coordinate $j$. The divisibility map $C_{N_k}\to C_{N_j}$ preserves edges and colors, so $f$ is a color-preserving graph homomorphism.

Conversely, suppose that a color-preserving graph homomorphism $f:H_{p_i}\to H_{p_j}$ exists. The component of $H_{p_i}$ with first coordinate $i$ is connected, so its image is contained in a component of $H_{p_j}$ with some first coordinate $k$. Then $k\le j$ and $p_k\le p_j$, while the induced color-preserving graph homomorphism $C_{N_i}\to C_{N_k}$ implies $N_k\mid N_i$. Hence $p_i\le p_k\le p_j$, proving \eqref{eq:hom-order}.
\end{proof}

The following lemma records the graph-theoretic structure required in the finite metric construction.

\begin{lemma}\label{lem:graph-hierarchy}
Let $Y=\{x_i:i\in\NN\}$ contain at least two points, and let $(\mathcal P_n)_{n\ge0}$ be a sequence of finite partitions of $Y$. Assume that $\mathcal P_0=\{Y\}$, that $\mathcal P_{n+1}$ refines $\mathcal P_n$, and that $\{x_i\}\in\mathcal P_n$ whenever $i\le n$. Put
\[
 \mathscr A=\bigcup_{n\ge0}\mathcal P_n\cup\{Y\setminus\{x_i\}:i\in\NN\},
\]
ordered by inclusion, and for $A\in\mathscr A$ let $H_A$ be the graph supplied by \cref{lem:poset-graphs}. Then one can choose properly 3-colored finite graphs $F_n$ and color-preserving graph homomorphisms $\Phi_k:F_{k+1}\to F_k$ $(k\ge0)$ with the following properties.

For every $n$,
\begin{align*}
 V(F_n)&=\{(A,v):A\in\mathcal P_n,\ v\in V(H_A)\},\\
 E(F_n)&=\bigl\{\{(A,u),(A,v)\}:A\in\mathcal P_n,\ \{u,v\}\in E(H_A)\bigr\},
\end{align*}
and $(A,v)$ has the color of $v$ in $H_A$. For $0\le k\le n$, set $\Phi_{n,n}=\operatorname{id}_{F_n}$ and $\Phi_{k,n}=\Phi_k\circ\cdots\circ\Phi_{n-1}$ when $k<n$. If $s\in V(F_n)$, write $s_k=\Phi_{k,n}(s)=(A_k(s),v_k)$. Then
\begin{equation}\label{eq:hierarchy-nesting}
 A_n(s)\subseteq A_{n-1}(s)\subseteq\cdots\subseteq A_0(s)=Y.
\end{equation}
Moreover, for $s,t\in V(F_n)$, $s_k=t_k$ implies $s_\ell=t_\ell$ for every $0\le\ell\le k$.
\end{lemma}

\begin{proof}
For $B\in\mathcal P_{k+1}$, there is a unique $A\in\mathcal P_k$ with $B\subseteq A$. Denote it by $\operatorname{par}_k(B)$. By \eqref{eq:hom-order},
$\operatorname{Hom}_{\mathrm{col}}(H_B,H_{\operatorname{par}_k(B)})\ne\varnothing$, so choose a color-preserving graph homomorphism $\varphi_{k,B}:H_B\to H_{\operatorname{par}_k(B)}$.

Define $F_n$ by the formulas in the statement, and for $B\in\mathcal P_{k+1}$ and $v\in V(H_B)$ set
\[
 \Phi_k(B,v)=\bigl(\operatorname{par}_k(B),\varphi_{k,B}(v)\bigr).
\]
Since each $\varphi_{k,B}$ preserves edges and colors, so does $\Phi_k$. The inclusions in \eqref{eq:hierarchy-nesting} follow from the definition of $\operatorname{par}_k$. If $s_k=t_k$, then $\Phi_{\ell,k}(s_k)=\Phi_{\ell,k}(t_k)$ for every $0\le\ell\le k$, and hence $s_\ell=t_\ell$.
\end{proof}

\begin{proposition}\label{prop:F-to-embedding}
If a separable bounded metric space $X$ satisfies condition $(F)$, then $X$ is isometrically embeddable into $(\M,\dGH)$.
\end{proposition}

\begin{proof}
The one-point case is trivial, so let $D=\diam X>0$. Condition $(F)$ is invariant under rescaling, since the threshold $\eta$ for $D^{-1}d$ corresponds to $D\eta$ for $d$. By the scaling identity in \cref{prop:gh-basic}, we may therefore assume $\diam X=1$.

Choose a dense sequence $(x_i)_{i\in\NN}$, allowing repetitions, and put $Y=\{x_i:i\in\NN\}$. Then $Y$ contains at least two points. Moreover, $\Gamma_{o,\eta}(Y)$ is an induced subgraph of $\Gamma_{o,\eta}(X)$, so $Y$ also satisfies condition $(F)$. Apply \cref{lem:partitions} to $Y$ and fix the resulting sequence $(\mathcal P_n)$ of partitions. We also retain the notation $H_x(A)$ from that lemma. Next apply \cref{lem:graph-hierarchy} and fix the finite colored graphs $F_n$, the maps $\Phi_{k,n}$, and the sets $A_k(s)$.

For $x\in Y$ and $n\ge1$, define the following quantities on $V(F_n)$. For differently colored vertices $s,t\in V(F_n)$, set
\[
 I_n(s,t)=\{k\in\{0,\ldots,n\}:\{s_k,t_k\}\in E(F_k)\},
\]
and define
\begin{equation}\label{eq:m-finite}
 m_{x,n}(s,t)=
 \max\bigl(\{d(x,A_k(s)):k\in I_n(s,t)\}\cup\{0\}\bigr).
\end{equation}
If $k\in I_n(s,t)$, then the definition of the edges of $F_k$ gives some $A\in\mathcal P_k$ and $u,v\in V(H_A)$ with $s_k=(A,u)$ and $t_k=(A,v)$. In particular, $A_k(s)=A_k(t)=A$.

For distinct vertices $s,t$ of the same color, $s_n=s\ne t=t_n$, so there is a least $j$ with $s_j\ne t_j$. Set
$j(s,t)=\min\{j\in\{0,\ldots,n\}:s_j\ne t_j\}$ and define
\begin{equation}\label{eq:delta-finite}
 \delta_{x,n}(s,t)=
 \begin{cases}
 1,&j(s,t)=0,\\
 \max\{H_x(A_{j(s,t)-1}(s)),2^{-j(s,t)-1}\},&j(s,t)\ge1.
 \end{cases}
\end{equation}
Set $\delta_{x,n}(s,s)=0$. Since $j(s,t)=j(t,s)$ and, when $j(s,t)\ge1$, one has $s_{j(s,t)-1}=t_{j(s,t)-1}$, the function $\delta_{x,n}$ is symmetric. The function $m_{x,n}$ is symmetric as well. Since $\diam Y\le1$, one has $0\le m_{x,n}(s,t)\le1$ for differently colored pairs.

Let $\gamma_{01}=\gamma_{10}=3$, $\gamma_{02}=\gamma_{20}=5$, and $\gamma_{12}=\gamma_{21}=7$. Define $p_{x,n}:V(F_n)\times V(F_n)\to[0,\infty)$ by
\begin{equation}\label{eq:p-finite}
 p_{x,n}(s,t)=
 \begin{cases}
 \delta_{x,n}(s,t),&s,t\text{ have the same color},\\
 \gamma_{ab}+m_{x,n}(s,t),&s,t\text{ have colors }a\ne b.
 \end{cases}
\end{equation}

\medskip\noindent\textit{Claim 1.}\quad{}For every $x\in Y$ and $n\ge1$, the function $p_{x,n}$ is a metric.

We first prove that $\delta_{x,n}$ satisfies the ultrametric inequality on each color class. Let $s,t,u$ be distinct vertices of the same color. Interchanging $s$ and $t$ if necessary, assume $q=j(s,u)\le j(t,u)$. If $q=0$, then $\delta_{x,n}(s,u)=1$, while \cref{lem:partitions}\textup{(a)} and \eqref{eq:delta-finite} give $\delta_{x,n}(s,t)\le1$. If $q\ge1$, then $s_\ell=u_\ell=t_\ell$ for $0\le\ell<q$, so $j(s,t)\ge q$. By \eqref{eq:hierarchy-nesting},
$A_{j(s,t)-1}(s)\subseteq A_{q-1}(s)$, and \cref{lem:partitions}\textup{(a)} together with \eqref{eq:delta-finite} yields
$\delta_{x,n}(s,t)\le\delta_{x,n}(s,u)$. Cases in which two points coincide are immediate.

Next, if $s,t$ have the same color and $v$ has a different color, then
\begin{equation}\label{eq:m-delta}
 |m_{x,n}(s,v)-m_{x,n}(t,v)|\le\delta_{x,n}(s,t).
\end{equation}
If $m_{x,n}(s,v)\ne m_{x,n}(t,v)$, interchange $s$ and $t$ if necessary so that
$m_{x,n}(s,v)>m_{x,n}(t,v)$. If $j(s,t)=0$, then \eqref{eq:m-delta} follows from $0\le m_{x,n}\le1$. Assume $j(s,t)\ge1$, and put $k=j(s,t)-1$ and $A=A_k(s)=A_k(t)$. Since $I_n(s,v)$ is finite, some $\ell\in I_n(s,v)$ realizes the maximum in \eqref{eq:m-finite}, with
$m_{x,n}(s,v)=d(x,A_\ell(s))>0$. If $\ell\le k$, then $s_\ell=t_\ell$, so $\ell\in I_n(t,v)$ and $A_\ell(s)=A_\ell(t)$, which would imply
$m_{x,n}(t,v)\ge m_{x,n}(s,v)$, a contradiction. Hence $\ell>k$. Since $\{s_\ell,v_\ell\}\in E(F_\ell)$ and $\Phi_{k,\ell}$ is a graph homomorphism,
$\{s_k,v_k\}=\{t_k,v_k\}\in E(F_k)$, so $m_{x,n}(t,v)\ge d(x,A)$. Moreover, $A_\ell(s)\subseteq A$, and \cref{lem:partitions}\textup{(c)} gives
$d(x,A_\ell(s))\le d(x,A)+H_x(A)$. Therefore
\[
 m_{x,n}(s,v)-m_{x,n}(t,v)
 \le H_x(A)\le\delta_{x,n}(s,t),
\]
which proves \eqref{eq:m-delta}.

Positive definiteness of $p_{x,n}$ follows from \eqref{eq:delta-finite} and $\gamma_{ab}\ge3$. For the triangle inequality, use the ultrametric inequality when all three vertices have the same color. When exactly two vertices have the same color, use \eqref{eq:m-delta} together with $p_{x,n}(s,t)\le1$. When all three colors occur, the three side lengths belong to $[3,4]$, $[5,6]$, and $[7,8]$, respectively. Hence the largest side is at most $8$, while the sum of the other two is at least $8$. This proves Claim 1.

For every $n$ and $1\le i\le n$, set
\[
 A_i^{(n)}=(V(F_n),2p_{x_i,n}).
\]

\medskip\noindent\textit{Claim 2.}\quad{}For every $1\le i,j\le n$,
\begin{equation}\label{eq:model-exact}
 \dGH(A_i^{(n)},A_j^{(n)})=d(x_i,x_j).
\end{equation}

By \cref{lem:partitions}\textup{(b)} and the $1$-Lipschitz property of $x\mapsto d(x,A)$, the functions $\delta_{x,n}(s,t)$ and $m_{x,n}(s,t)$ are $1$-Lipschitz in $x$ on their respective domains. Hence
$|p_{x,n}(s,t)-p_{y,n}(s,t)|\le d(x,y)$. The identity correspondence therefore has distortion at most $2d(x_i,x_j)$, and \eqref{eq:gh-corr} gives
$\dGH(A_i^{(n)},A_j^{(n)})\le d(x_i,x_j)$.

For the reverse inequality, set $\Delta=d(x_i,x_j)$. If $\Delta=0$, equality follows from the upper bound. Assume $\Delta>0$. Suppose that there exists a correspondence
$R\in\operatorname{Corr}(A_i^{(n)},A_j^{(n)})$ with $\dis R<2\Delta$. If $(s,s'),(t,t')\in R$ and $s,t$ have the same color, then $s',t'$ also have the same color. Indeed, otherwise
$2p_{x_i,n}(s,t)\le2$ and $2p_{x_j,n}(s',t')\ge6$, so $\dis R\ge4>2\Delta$. Thus the color of a vertex corresponding to a vertex of color $a$ depends only on $a$. Applying the same argument to $R^{-1}$ shows that this defines a permutation $\sigma$ of $\{0,1,2\}$.

If $(s,s'),(t,t')\in R$ and $s,t$ have different colors, then the distortion bound gives $|p_{x_i,n}(s,t)-p_{x_j,n}(s',t')|<\Delta\le1$. The three possible distance intervals $[3,4]$, $[5,6]$, and $[7,8]$ are separated by at least $1$, so
$\gamma_{ab}=\gamma_{\sigma(a)\sigma(b)}$. Since $\gamma_{01}$, $\gamma_{02}$, and $\gamma_{12}$ are distinct, $\sigma$ is the identity. Therefore every $(s,s')\in R$ consists of vertices of the same color.

Since $n\ge j$, the singleton $\{x_j\}$ belongs to $\mathcal P_n$. Let $G_{x_j}$ be the subgraph of $F_n$ with vertex set $\{\{x_j\}\}\times V(H_{\{x_j\}})$ and edge set
\[
 E(G_{x_j})=\bigl\{\{(\{x_j\},u),(\{x_j\},v)\}:\{u,v\}\in E(H_{\{x_j\}})\bigr\}.
\]
The map $\iota_{x_j}:V(H_{\{x_j\}})\to V(G_{x_j})$, $\iota_{x_j}(u)=(\{x_j\},u)$, is a color-preserving graph isomorphism.

If $\{s,t\}\in E(G_{x_j})$, then $\Phi_{k,n}$ is a graph homomorphism, so $\{s_k,t_k\}\in E(F_k)$ for every $0\le k\le n$. Moreover, $A_n(s)=\{x_j\}\subseteq A_k(s)$, so $d(x_i,A_k(s))\le\Delta$, with equality when $k=n$. Hence $m_{x_i,n}(s,t)=\Delta$.

For each $s\in V(G_{x_j})$, choose $q(s)\in V(F_n)$ with $(s,q(s))\in R$. The map $q$ preserves colors. For every edge $e=\{s,t\}\in E(G_{x_j})$,
\[
 |\Delta-m_{x_j,n}(q(s),q(t))|
 =|p_{x_i,n}(s,t)-p_{x_j,n}(q(s),q(t))|<\Delta,
\]
so $m_{x_j,n}(q(s),q(t))>0$. Define
\[
 J_e=\{r\in\{0,\ldots,n\}:\{q(s)_r,q(t)_r\}\in E(F_r),\ d(x_j,A_r(q(s)))>0\}.
\]
By \eqref{eq:m-finite}, the set $J_e$ is nonempty.

Let $C$ be a connected component of $G_{x_j}$ containing an edge, and put
$k_C=\min\bigcup_{e\in E(C)}J_e$. For any edge $e=\{s,t\}\in E(C)$, choose $r\in J_e$. Since $k_C\le r$ and $\Phi_{k_C,r}$ sends edges to edges,
$\{q(s)_{k_C},q(t)_{k_C}\}\in E(F_{k_C})$. Thus $s\mapsto q(s)_{k_C}$ is a color-preserving graph homomorphism from $C$ to $F_{k_C}$. Edges of $F_{k_C}$ never join vertices with different first coordinates. Since $C$ is connected, there exists $A_C\in\mathcal P_{k_C}$ such that
$q(s)_{k_C}=(A_C,v_s)$ for every $s\in V(C)$. By the definition of $k_C$, some edge $e=\{s,t\}$ satisfies $k_C\in J_e$, and hence $d(x_j,A_C)>0$. Therefore $A_C\subseteq Y\setminus\{x_j\}$. By \eqref{eq:hom-order}, there is a color-preserving graph homomorphism $H_{A_C}\to H_{Y\setminus\{x_j\}}$. Composing it with $s\mapsto v_s$ gives a color-preserving graph homomorphism $C\to H_{Y\setminus\{x_j\}}$.

If a connected component $C$ consists of a single isolated vertex $s$, choose a vertex of $H_{Y\setminus\{x_j\}}$ having the same color as $s$. Such a vertex exists because $H_{Y\setminus\{x_j\}}$ uses all three colors. For each connected component $C$ of $G_{x_j}$, fix a color-preserving graph homomorphism
$\psi_C:C\to H_{Y\setminus\{x_j\}}$ as constructed above. Define $\Psi(s)=\psi_C(s)$ for $s\in V(C)$. Since the connected components partition $V(G_{x_j})$, this defines a color-preserving graph homomorphism
$\Psi:G_{x_j}\to H_{Y\setminus\{x_j\}}$. Consequently,
$\Psi\circ\iota_{x_j}:H_{\{x_j\}}\to H_{Y\setminus\{x_j\}}$ is also a color-preserving graph homomorphism. By \eqref{eq:hom-order}, this would imply $\{x_j\}\subseteq Y\setminus\{x_j\}$, a contradiction. Hence no correspondence of distortion less than $2\Delta$ exists, and
$\dGH(A_i^{(n)},A_j^{(n)})\ge\Delta$. Claim 2 is proved.

\medskip\noindent\textit{Claim 3.}\quad{}For each fixed $i$, the family $(A_i^{(n)})_{n\ge i}$ is uniformly totally bounded.

For differently colored pairs, $0\le m_{x,n}\le1$. For pairs of the same color, \cref{lem:partitions}\textup{(a)} and \eqref{eq:delta-finite} give $0\le\delta_{x,n}\le1$. Hence $p_{x,n}\le8$, and therefore $\diam A_i^{(n)}\le16$. Fix $i$ and $r>0$, and choose an integer $k\ge i$ such that $2^{1-k}<r$. Let $n\ge k$. For each $w\in V(F_k)$, set
$C_w=\Phi_{k,n}^{-1}(\{w\})$. If distinct $s,t\in C_w$, then $s_k=t_k$, so $s,t$ have the same color and $j(s,t)\ge k+1$. By \eqref{eq:hierarchy-nesting}, \eqref{eq:partition-decay}, and $k\ge i$,
\[
 H_{x_i}(A_{j(s,t)-1}(s))\le H_{x_i}(A_k(s))\le2^{-k},
 \qquad 2^{-j(s,t)-1}\le2^{-k-2}.
\]
Thus
$d_{A_i^{(n)}}(s,t)=2\delta_{x_i,n}(s,t)\le2^{1-k}<r$. Choosing one point from each nonempty $C_w$ gives an $r$-net with at most $|V(F_k)|$ points. Including the finitely many initial terms $i\le n<k$, the cardinalities of the required $r$-nets for all $n\ge i$ are bounded by
$\max_{i\le m\le k}|V(F_m)|$. This proves Claim 3.

By Claims 1--3, the finite metric spaces $A_i^{(n)}$ satisfy the hypotheses of \cref{lem:finite-approx}. Hence the normalized space $X$ is isometrically embeddable into $\M$. Multiplying the metrics on all compact representatives by $D$ yields an isometric embedding of $X$ equipped with its original metric.
\end{proof}

\begin{proof}[Proof of \cref{thm:main}]
By \cref{prop:intrinsic}, conditions \textup{(ii)}--\textup{(iv)} are equivalent. Proposition \ref{prop:embedding-to-ball} gives \textup{(i)}$\Rightarrow$\textup{(iii)}. Proposition \ref{prop:F-to-embedding} gives \textup{(ii)}$\Rightarrow$\textup{(i)}. Hence all four conditions are equivalent.
\end{proof}

The bounded case is therefore complete. The corresponding question without boundedness remains open.

\begin{problem}\label{prob:unbounded}
Characterize the metric spaces that are isometrically embeddable into $(\M,\dGH)$ without assuming boundedness.
\end{problem}

\section*{Acknowledgments}
The author would like to thank OpenAI for providing ChatGPT for Academic Researchers. The author used the GPT-6.0 and GPT-6.2 models to assist in preparing this manuscript. The author reviewed and revised the mathematical content and takes full responsibility for the final manuscript.

\bigskip
\noindent\textsc{Mathematical Institute, Tohoku University}\\
Sendai 980-8578, Japan\\
\textit{Email address:} \texttt{fukano.ryo.t5@dc.tohoku.ac.jp}

\end{document}